\documentclass{amsart}
\usepackage[british]{babel}
\usepackage{amssymb,amsfonts}
\usepackage{graphicx,enumerate,xcolor}
\usepackage{scalerel}
\usepackage{stackengine,wasysym}
\usepackage{geometry}
\catcode`\<=\active \def<{
\fontencoding{T1}\selectfont\symbol{60}\fontencoding{\encodingdefault}}
\catcode`\>=\active \def>{
\fontencoding{T1}\selectfont\symbol{62}\fontencoding{\encodingdefault}}
\newcommand{\assign}{:=}
\newcommand{\longleftarrowlim}{\mathop{\longleftarrow}\limits}
\newcommand{\longrightarrowlim}{\mathop{\longrightarrow}\limits}

\newcommand{\rightarrowlim}{\mathop{\rightarrow}\limits}
\newcommand{\nocomma}{}
\newcommand{\nosymbol}{}

\newcommand{\tmdummy}{$\mbox{}$}
\newcommand{\tmem}[1]{{\em #1\/}}
\newcommand{\tmmathbf}[1]{\ensuremath{\boldsymbol{#1}}}
\newcommand{\tmop}[1]{\ensuremath{\operatorname{#1}}}
\newcommand{\tmstrong}[1]{\textbf{#1}}

\newenvironment{enumeratenumeric}{\begin{enumerate}[1.] }{\end{enumerate}}
\newenvironment{enumerateroman}{\begin{enumerate}[i.] }{\end{enumerate}}
\newenvironment{itemizeminus}{\begin{itemize} }{\end{itemize}}

\theoremstyle{plain}
\newtheorem{theorem}{Theorem}
\newtheorem{proposition}{Proposition}
\newtheorem{corollary}{Corollary}
\newtheorem{lemma}{Lemma}

\newtheorem{definition}{Definition}

\theoremstyle{definition}
\newtheorem{example}{Example}
\newtheorem{remark}{Remark}
\newtheorem*{notation}{Notation}

\newcommand{\eff}{\ensuremath{\mathcal{E}f \nocomma f}}

\newcommand{\eeq}[2]{| #1 {\approx} #2 |}

\newcommand{\app}[2]{\ensuremath{#1 \cdot #2}}

\newcommand{\ent}{\: \vdash \:}
\newcommand{\pn}{\mathcal{P}\mathbb{N}}
\newcommand{\isval}[1]{\pn {\vDash} #1}
\newcommand{\eset}[1]{\left(#1, {\approx}\right)}

\newcommand{\rep}[1]{\ensuremath{\mathfrak{R}}_{#1}}

\newcommand{\pao}[1]{\ensuremath{#1^{\langle I \rangle}}}

\newcommand{\nno}{\ensuremath\tmop{{\mathbf{N}}}}
\newcommand{\htpy}{{\looparrowright}}

\newcommand{\astt}{{\star}}

\newcommand{\rev}{\ensuremath{\tmop{\mathbf{rev}}}}
\newcommand{\pad}{\ensuremath{\tmop{\mathbf{pad}}}}
\newcommand{\map}{\ensuremath{\tmop{\mathbf{map}}}}
\newcommand{\pco}{\ensuremath{\tmop{Path}}}
\newcommand{\flatt}[2]{#1^{\left(#2\right)}_{\flat}}
\newcommand{\pic}[1]{\left\langle#1\right\rangle}

\newcommand{\src}{ \ensuremath{\partial^-}}
\newcommand{\tgt}{\ensuremath{\partial^+}}
\newcommand{\piq}[2]{\left\langle#1, #2\right\rangle}

\newcommand{\deq}{\overset{def.}{=}}

\newcommand{\idle}[1]{{\star}}

\newcommand{\llift}[1]{\ensuremath{^{\pitchfork} #1}}
\newcommand{\rlift}[1]{\ensuremath{#1^{\pitchfork}}}

\newcommand{\llp}[1]{\ensuremath{^{\pitchfork} #1}}
\newcommand{\rlp}[1]{\ensuremath{#1^{\pitchfork}}}

\newcommand\reallywidetilde[1]{\ThisStyle{%
  \setbox0=\hbox{$\SavedStyle#1$}%
  \stackengine{-.1\LMpt}{$\SavedStyle#1$}{%
    \stretchto{\scaleto{\SavedStyle\mkern.2mu\AC}{.5150\wd0}}{.6\ht0}%
  }{O}{c}{F}{T}{S}%
}}

\begin{document}

\title{Hurewicz fibrations in elementary toposes}

\author{Krzysztof Worytkiewicz}

\date{August 30, 2019}

\thanks{
This work is part of the output of the AmSud Math project
\tmem{Logic, Categories and Complexity} and as such benefitted from the
relevant financial support. This author wishes to thank the people responsible
for coordination at AmSud Math for their efficiency and kindness. He
acknowledges and thanks Benno van den Berg, F{\'e}lix Castro, Martin Hyland,
Jan Mina{\v c}, Alexandre Miquel and Jaap van Oosten.
}

\begin{abstract}
We study formal counterparts of Hurewicz fibrations and related topological
notions in elementary toposes with NNO. The constructions are based on a
specific notion of interval and lead to a structure of category of fibrant
objects on toposes equipped with such a datum. We get in fact slightly more
as the building blocks are derived from a weak factorisation system.
\end{abstract}

{\maketitle}

\section{Introduction}
The fundamental geometric notion of {\tmem{cohe{\tmem{}}sion}}, a way of
distinguishing parts of a whole, has been traced back in
{\cite{bell2009cohesiveness}} to Aristotle's {\tmem{Categories, Book VI}}. It
is a recurrent theme in the present work, which originates in a specific
notion of cohesion encountered in realisability toposes over partial
combinatory algebras: two points {\tmem{stick}} if their sets of realisers
intersect non-trivially {\cite{van2015notion}}. Unfortunately, the established
theory of {\tmem{axiomatic cohesion}}
{\cite{lawvere1994cohesive,lawvere2007axiomatic,marmolejo2017relation}} does
not cover this case, the relevant adjoints being the other way around. Given
our intended applications, rather than attempting to find a notion of cohesion
suitable for realisability toposes, we choose to encompass it into a notion of
{\tmem{interval}} acting as a bootstrapping cohesion datum. In particular, our
notion of interval yields a simplicial resolution of any object. This
simplicial resolution, a nerve construction which we call {\tmem{path
complex}}, encodes cohesion in the sense that its inhabiting paths stay in the
same connected component {\tmem{by construction}}. In a realisability topos
equipped a bipointed object such that its points stick (known as $\Delta
\tmmathbf{2}$), this means that paths can only be laid out along points that
stick. Connected components are thus identical to path-connected ones in this
setting. A weak notion of geometric realisation turns such a path complex into
a {\tmem{path object}}. Quite surprisingly, this is enough to have workable
internal versions of constructs known from topology: Hurewicz fibrations,
fundamental categories, homotopy and strong deformation retracts among others.
These are the ingredients giving rise to a weak factorisation system and
further to a structure of category of fibrant objects
{\cite{brown1973abstract}} on a topos with NNO equipped with such an interval.
This author likes to think of this material as a generalisation of van
Oosten's work {\cite{van2015notion}} on the effective topos
{\cite{hyland1982effective}}.

Unsurprisingly, our leading example will be latter. In order not to overload
the exposition, some of the relevant details are to be found in an appendix.
In Section \ref{sec:toposes} we recall what a topos is and briefly review some
relevant features. In Section \ref{sec:path-complex} we introduce an
elementary notion of interval in a topos and successively add features
required to build the path complex. In Section \ref{sec:path-object} we
introduce a weak notion of geometric realisation of a path complex, which
consists of just modding out degeneracies without gluing along adjacent faces.
It turns out that the result of the construction yields a functorial notion of
an internal category: the {\tmem{fundamental category}} of an object. In this
context, the fundamental category plays the r{\^o}le of a {\tmem{path
object}}: we use it to define (right) homotopy and homotopy equivalence. We
also describle how {\tmem{contracting homotopy}}
{\cite{jardine2006categorical,van2012topological}} arises in this context. At
this point we also introduce the {\tmem{Hurewicz property}} of an interval. It
is a technical condition of homotopical nature on the associated path complex.
In Section \ref{sec:hurewicz-fibs} we introduce the central notion of
{\tmem{Hurewicz fibration}}, which turns out to formally behave like it's
topological counterpart, and also present some relevant instances. In Section
\ref{sec:wfs}, in a sense the technical crux of the paper, we introduce
{\tmem{strong deformation insertions}}, that is insertions of strong
deformation retracts, and show that the class of the latter and the class of
Hurewicz fibrations form a weak factorisation system. In Section \ref{sec:cfo}
we fill in the remaining gaps in order to exhibit the structure of
{\tmem{category of fibrant objects}} {\cite{brown1973abstract}} on a topos
equipped with a Hurewicz interval.

We have in fact a bit more than Brown's original formulation here, given that
the class of Hurewicz fibrations is part of a weak factorisation system. It
would be tempting to conjecture that this is in fact one half of a ``Str{\o}m
model structure'' to be unveiled. On the other hand, the present structure of
category of fibrant objects makes constructions like universal bundles and
more generally cohomological techniques available. An interesting application
of this cercle of ideas would be the construction of models of Homotopy Type
Theory in toposes equipped with a Hurewicz interval, this since categorical
models of HoTT which are {\tmem{tribes}} {\cite{joyal2017notes}} and
categories of fibrant objects are in a tight relationship given by a
DK-equivalence {\cite{kapulkin2019internal}}. This would in particular give
rise to realizability models of HoTT.

\section{Preliminaries\label{sec:toposes}}
\begin{notation}
Given a category $\mathbb{C}$ we shall write $\mathbb{C}_0$ for
its class of objects and $\mathbb{C}_1$ for its class of morhisms. $X \in
\mathbb{C}$ stands for $X \in \mathbb{C}_0$. We shall write
\begin{itemizeminus}
  \item  $X' \vartriangleleft X$ to indicate that $X$ is a subobject of $X$;

  \item $\star$ for a terminal object.
\end{itemizeminus}
Given a class of morhisms $\mathcal{A} \subset \mathbb{C}_1$ we shall write
$\llift{\mathcal{A}}$ for it's left lifting class and $\rlift{\mathcal{A}}$ for
its right lifting class.
\end{notation}

\begin{definition}
  {\tmdummy}

  \begin{enumeratenumeric}
    \item Assume a category $\mathbb{C}$ with finite products. $\mathbb{C}$ is
    {\tmem{cartesian closed}} if the functor $(-) \times Y : \mathbb{C}
    \rightarrow \mathbb{C}$ has a right adjoint $(-)^Y : \mathbb{C}
    \rightarrow \mathbb{C}$ for all $Y \in \mathbb{C}$. Evaluating the latter
    is called {\tmem{exponentiation}}.

    \item Assume a category $\mathbb{C}$ with finite limits. A
    {\tmem{subobject classifier}} in $\mathbb{C}$ is a mono $\tmop{tt} : \star
    \rightarrow \Omega$ from the terminal object $\star$ such that for any
    mono $a : A \rightarrowtail X$ there is a unique {\tmem{classifying
    morphism}} $\chi_A : X \rightarrow \Omega$ such that there is a pullback
    diagram

    \begin{center}
      \includegraphics[trim = {38 0 0 30}]{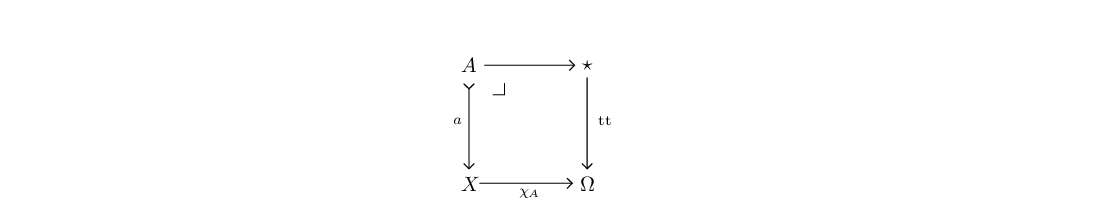}
    \end{center}

    \item A {\tmem{topos}} is a finitely complete $\tmop{CCC}$ with a
    subobject classifier.
  \end{enumeratenumeric}
\end{definition}

{\noindent}This low-key definition entails a vaste array of consequences which
have been extensively studied since the 1970's
{\cite{kock1971elementary,johnstone1974associated,lawvere1976variable}}. In
particular, constructions in a topos $\mathbb{T}$ can be performed in many
relevant cases using $\mathbb{T}$'s {\tmem{internal logic}} built on top of
$\mathbb{T}$'s {\tmem{internal language}}
{\cite{osius1974internal,fourman1977logic,boileau1981logique,maclane1975sets}},
that is using a type theory formally akin to a (in general constructive)
calculus of sets and functions. It is for instance the case that given an
object $X \in \mathbb{T}$ \ the exponential $\Omega^X$ formally behaves like a
{\tmem{powerset}}, whence the name {\tmem{$X$'s powerobject}}. \ We refer the
reader to
{\cite{johnstone2002sketches,maclane2012sheaves,borceux1994handbook}} for the
lore. Notice that what we choose to call {\tmem{topos}} here often goes under
{\tmem{elementary topos}} in the litterature.

\begin{definition}
  Assume a topos $\mathbb{T}$. A {\tmem{natural numbers object}} $\nno \in
  \mathbb{T}$ is part of the sequence $\star \longrightarrowlim^Z \nno
  \longrightarrowlim^S \nno$ which is an initial object in the category of
  sequences $\star \rightarrow X \rightarrow X$ in $\mathbb{T}$.
\end{definition}

{\noindent}Assume $\mathbb{T}$ a topos with NNO $\nno$.

\begin{remark}
  Assume a morphism $\phi : X \rightarrow \Omega^Y$. The subobject of dependent
  pairs
  \begin{eqnarray*}
    \bigsqcup_{x : X} \phi (x) & \deq & \{ (x, y) : X \times Y|y \in \phi (x)
    \}
  \end{eqnarray*}
  is classified by $\phi$'s exponential transpose $\check{\phi
  :} X \times Y \rightarrow \Omega$. It reflects the set-theoretical intuition
  of a relation $X \nrightarrow Y$ as a disjoint union of a family of subsets
  of $Y$ indexed by $X$.
\end{remark}

\begin{definition}
  Let $X, Y \in \mathbb{T}$. A family of subobjects of $Y$ indexed by $X$ is a
  morphism $\phi : X \rightarrow \Omega^Y$.
\end{definition}

\begin{notation}
We shall write $(Y_x)_{x : X}$ for a family $\phi : X \rightarrow
\Omega^Y$ if $\phi$ is understood, and \ $\bigsqcup_{x : X} Y_x$
accordingly.
\end{notation}

\begin{definition}
  Let $X \in \mathbb{T}$. The object
  \begin{eqnarray*}
    \tmop{List} (X) & \deq & \left\{ u : \left( \astt + X \right)^{\nno} |
    \exists n : \nno \nosymbol . \forall k : \nno . (k < n \Rightarrow u_k \in
    X) \wedge (k \geqslant n \Rightarrow u_k = \star) \right\}
  \end{eqnarray*}
  is called {\tmem{list object over }}$X$.
\end{definition}

\begin{remark}
  {\tmdummy}

  \begin{enumeratenumeric}
    \item A topos being a CCC, there is the morphism $\map : Y^X \rightarrow
    \tmop{List} (X) \rightarrow \tmop{List} (Y)$ for all $X, Y \in
    \mathbb{T}$.

    \item An NNO being decidable, there is a length morphism $\ell :
    \tmop{List} (X) \rightarrow \nno$ for all $X \in \mathbb{T}$ and a list
    reversing isomorphism $\wp : \tmop{List} (X) \longrightarrowlim^{\cong}
    \tmop{List} (X)$.
  \end{enumeratenumeric}
\end{remark}

\section{The path complex}\label{sec:path-complex}

\begin{definition}
  An object $X \in \mathbb{C}$ in a category $\mathbb{C}$ is
  {\tmem{well-pointed}} if, given (arbitrary) morphisms $f, g : X \rightarrow
  Y$, $f (x) = g (x)$ for all global sections $x : \astt \rightarrow X$
  implies $f = g$.
\end{definition}

\begin{remark}
  A boolean topos can be characterised as a topos where every object is
  well-pointed.
\end{remark}

\begin{definition}
  $I \in \mathbb{T}$ is an {\tmem{elementary interval}} provided it
  \begin{enumerateroman}
    \item is well-pointed;

    \item has precisely two global sections $\#0, \#1 : \astt \rightarrow I$.
  \end{enumerateroman}
\end{definition}

\begin{remark} \label{example:I}
  In {\eff}, a global section $\astt
  \rightarrow (X, \approx)$ selects an equality class in $\underline{X} /
  \approx$ (c.f. Remark \ref{rem:Delta}). Let $I$ be the assembly
  \[ I \deq  (\{ i_0, i_1 \} ; E (i_0) = E (i_1) =\mathbb{N}) \]
  so $I \cong \Delta \tmmathbf{2}$. Assume a morphism $u : I \rightarrow (X,
  \approx)$. The total condition entails that there are elements $x_0, x_1 \in
  X$ such that there is a Turing machine uniformely realising
  \begin{eqnarray*}
    E (i_0) & \rightarrow & \mathfrak{R}_u (i_0, x_0) \\
    E (i_1) & \rightarrow & \mathfrak{R}_u (i_1, x_1)
  \end{eqnarray*}
  while the other conditions entail
  \begin{eqnarray*}
    \mathfrak{R}_u (i_0, x) & \leftrightarrow & | x \approx x_0 | \\
    \mathfrak{R}_u (i_1, x) & \leftrightarrow & | x \approx x_1 |
  \end{eqnarray*}
  Hence a morphism $I \rightarrow (X, \approx)$ determines and is
  determined by global sections $$\lceil x \rceil, \lceil x' \rceil : \astt
  \rightarrow (X, \approx)$$ verifying
  \begin{eqnarray*}
    \left( \bigcup_{y \approx x} E (y) \right) \cap \left( \bigcup_{y' \approx
    x'} E (y') \right) & \neq & \varnothing
  \end{eqnarray*}
(c.f. {\cite{hyland1982effective}} sec. 3). $I$ is in particular
  well-pointed, so it is an elementary interval.
\end{remark}

\begin{definition}
  Let $I$ be an elementary interval in $\mathbb{T}$. {\tmem{Elementary
  intervals of length $n$}} are obtained by gluing copies of $I$
  \begin{eqnarray*}
    I_0 & \deq & \astt\\
    I_{n + 1} & \deq & I_n +_1 I
  \end{eqnarray*}
  by pushout

  \begin{center}
    {\includegraphics[trim = {38 0 0 30}]{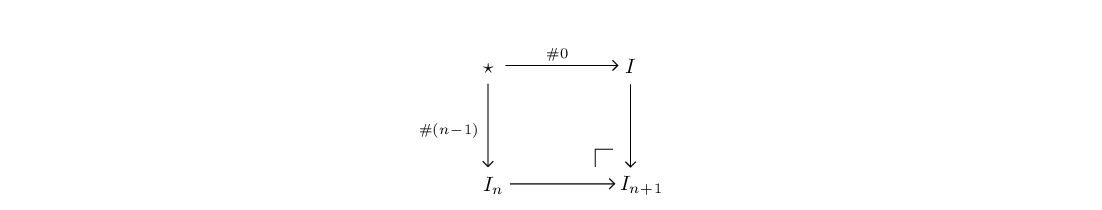}}
  \end{center}

  {\noindent}Let $X \in \mathbb{T}$.We shall call a morphism $I_n \rightarrow
  X$ {\tmem{path of degree $n$}} (in $X$) so $X^{I_n}$ is the {\tmem{object of
  paths of degree $n$}} (in $X$).
\end{definition}

\begin{remark}
  It can be shown by induction that $I_n$ is well-pointed.
\end{remark}

\begin{remark}
  {\label{rem:rpath}}
  In {\eff} we can
  construct $I_n$ as the assembly
  \[ (\{ i^{(n)}_0, \cdots, i^{(n)}_{n - 1} \} ; \forall 0 \leqslant k < n. E
     (i^{(n)}_k) =\mathbb{N}) \]
  so $I_n \cong \Delta \{ 0, \cdots, n - 1 \}$. Assume a morphism $u : I_n
  \rightarrow (X, \approx)$. Similarly to the case $n = 2$ (c.f. Remark
  \ref{example:I}), the total condition entails that there are elements $x_0,
  \cdots, x_{n - 1} \in X$ such that there is a Turing machine uniformely
  realising
  \begin{eqnarray*}
    E (k) & \rightarrow & \mathfrak{R}_u (i_k, x_k)
  \end{eqnarray*}
  for $k \in \{ 0, \cdots, n - 1 \}$, while the other conditions entail
  \begin{eqnarray*}
    \mathfrak{R}_u (i_k, x) & \leftrightarrow & | x \approx x_k |
  \end{eqnarray*}
  for $k \in \{ 0, \cdots, n - 1 \}$. A morphism $u : I_n \rightarrow (X,
  \approx)$ thus determines and is determined by a list $[\lceil x_0 \rceil ;
  \cdots, \lceil x_{n - 1} \rceil]$ of global sections of $(X, \approx)$ such
  that for any ordered subset $S \subset \{ 0, \cdots, n - 1 \}$
  \begin{eqnarray*}
    \bigcap_{s \in S} \left( \bigcup_{y \approx x_s} E (y) \right) & \neq &
    \varnothing
  \end{eqnarray*}
  Intuitively, such a morphism is a {\tmem{contractible path}}.
\end{remark}

\begin{remark}
  \label{rem:lists}Let $X \in \mathbb{T}$ and $n : \nno$. We have
  \begin{eqnarray*}
    X^{I_n} & \cong & \left\{ l : \tmop{List} (X) | (\exists \omega : X^{I_n}
    . \forall 0 \leqslant i < n. l (i) = \omega_{\#i}) \wedge \forall k : \nno
    . \left( x_k \in \star \text{ } \Leftrightarrow \text{ } k \geqslant n
    \right) \right\}
  \end{eqnarray*}
  since $I_n$ is well-pointed, so in particular $X^{I_n} \vartriangleleft
  \tmop{List} (X)$ for all $n : \nno$. We thus have a family $(X^{I_n})_{n :
  \nno}$ of subobjects of $\tmop{List} (X)$.
\end{remark}

\begin{notation}
When convenient, we shall use the list notation $[w_{\#0} ; \cdots
; w_{\# (n - 1)}]$ for a path of degree $n$.
\end{notation}

\begin{remark}
  Let $I$ be an elementary interval. For any $n \geqslant 0$ and $0 \leqslant
  i \leqslant n$ there is the {\tmem{$i$-th coface function}}
  \begin{eqnarray*}
    \delta^{(i)} : \Gamma (I_n) & \longrightarrow & \Gamma (I_{n + 1})\\
    \#j & \mapsto & \left\{\begin{array}{lll}
      \#j &  & j < i\\
      \# (j + 1) &  & j \geqslant i
    \end{array}\right.
  \end{eqnarray*}
  Similarly, for any $n \geqslant 1$ and $0 \leqslant i \leqslant n - 1$ there
  is the {\tmem{$i$-th codegeneracy function}}
  \begin{eqnarray*}
    \sigma^{(i)} : \Gamma (I_{n + 1}) & \longrightarrow & \Gamma (I_n)\\
    \#j & \mapsto & \left\{\begin{array}{lll}
      \#j &  & j \leqslant i\\
      \# (j - 1) &  & j > i
    \end{array}\right.
  \end{eqnarray*}
\end{remark}

\begin{definition}
  An elementary interval $I$ is {\tmem{cosimplical}} provided coface functions
  $\delta^{(i)}$ and codegeneracy functions $\sigma^{(i)}$ uniquely determine
  morphisms $\delta_i : I_n \rightarrow I_{n + 1}$ and $\sigma_i : I_{n + 1}
  \rightarrow I_n$. We shall call these morphisms {\tmem{elementary cofaces}}
  and {\tmem{elementary codegeneracies}}, respectively. Moreover, we shall
  call {\tmem{cofaces}} respectively {\tmem{codegeneracies}} compositions of
  the elementary ones.
\end{definition}

\begin{notation}
Assume a computable expression $t (x)$ with $x \in \tmop{FV} (t)$.
We shall write $\Lambda x.t$ for the code of the Turing machine corresponding
to the partial recursive function $\lambda x.t$.
\end{notation}

\begin{remark}
  \label{ex:cosimp}The elementary interval $\Delta \tmmathbf{2}$ in {\eff} is
  cosimplicial. The $I_n$'s are assemblies, so a global section $\star
  \rightarrow I_n$ is uniquely determined by an element of the underlying set
  while a morphism $f : I_m \rightarrow I_n$ is uniquely determined by a
  tracked function on the underlying sets. The {\tmem{$i$-th coface function}}
  \begin{eqnarray*}
    \delta^{(i)} : \Gamma (I_n) & \rightarrow & \Gamma (I_{n + 1})
  \end{eqnarray*}
  admits the tracker
  \[ \Lambda j. \text{{\tmstrong{if}}} \; j < i \; \text{{\tmstrong{then}}} \; j
     \; \text{{\tmstrong{else}}} \;  j + 1 \]
  when seen as a function $\{ i_0, \cdots, i_{n - 1} \} \rightarrow \{ i_0,
  \cdots, i_n \}$. Similarly, the {\tmem{$i$-th codegenacy function}}
  \begin{eqnarray*}
    \sigma^{(i)} : \Gamma (I_{n + 1}) & \longrightarrow & \Gamma (I_n)
  \end{eqnarray*}
  admits the tracker
  \[ \Lambda j. \text{{\tmstrong{if}}} \; j \leqslant i \; \text{{\tmstrong{then}}}
    \; j \; \text{{\tmstrong{else}}} \; j - 1 \]
  when seen as a function $\{ i_0, \cdots, i_n \} \rightarrow \{ i_0, \cdots,
  i_{n - 1} \}$.
\end{remark}

\begin{remark}
  Let $I$ be a cosimplicial interval and $\mathbb{I} \subset \mathbb{T}$ be
  the subcategory with objects the $I_n$'s and monotone morphisms, the latter
  are generated by elementary cofaces and codegeneracies modulo cosimplicial
  identities. $\mathbb{I}$ is a monoidal category with tensor given by pushout

  \begin{center}
    {\includegraphics[trim = {38 0 0 30}]{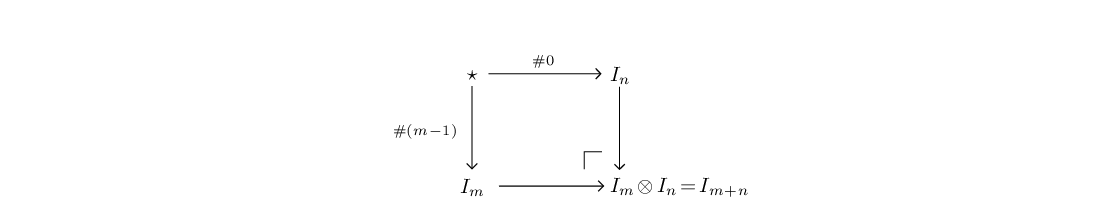}}
  \end{center}

  {\noindent}Any morphism in $\mathbb{I}$ admits a normal form. $\mathbb{I}$
  is in fact equivalent to $\Delta^+$, the augmented simplicial category.
\end{remark}

\begin{remark}
  \label{rem:eilenberg-zilber}Let $X \in \mathbb{T}$. The family $\pco (X)
  \deq (X^{I_n})_{n : \mathbb{N}}$ is a simplicial object with faces and
  degeneracies given by precomposition
  \begin{eqnarray*}
    d_i & = & \lambda w : \pco (X)_n . w \circ \delta_i : \tmop{Path}_n (X)
    \rightarrow \tmop{Path}_{n - 1} (X)\\
    s_i & = & \lambda w : \pco (X)_n . w \circ \sigma_i : \tmop{Path}_n (X)
    \rightarrow \tmop{Path}_{n + 1} (X)
  \end{eqnarray*}
  respectively. We shall call this simplicial object {\tmem{the path complex
  of $X$}} (with respect to $I$). Notice that $\pco (X)_0 = X^{I_0} \cong X$
  as $I_0 \cong \star$.
\end{remark}

\section{Path object and homotopy}\label{sec:path-object}

\begin{remark}
  Let $\underline{\tmop{Path}} (X)  \deq  \bigsqcup_{n : \mathbb{N}} \pco
  (X)_n$ so
  \begin{eqnarray*}
    \underline{\tmop{Path}} (X) & = & \left\{ (n, w) : \nno \times \tmop{List}
    (X) |w \in X^{I_n} \right\}
  \end{eqnarray*}
  (c.f. Remark \ref{rem:lists}). As $\nno$ is decidable, a face
  $d : \tmop{Path} (X)_m \rightarrow
  \tmop{Path} (X)_n$ induces an endomorphism
  \begin{eqnarray*}
    \underline{d} : \underline{\tmop{Path}} (X) & \rightarrow &
    \underline{\tmop{Path}} (X)
  \end{eqnarray*}
  constructed by the term
  $$
  \lambda (k, w) : \underline{\tmop{Path}} (X) .
  \text{{\tmstrong{if}}} \; k = m \; \text{{\tmstrong{then}}} \;
  (n,d (w)) \; \text{{\tmstrong{else}}} \; (k, w)
  $$
  Similarly, a degeneracy $s$ induces an endomorphism
  \begin{eqnarray*}
    \underline{s} : \underline{\tmop{Path}} (X) & \rightarrow &
    \underline{\tmop{Path}} (X)
  \end{eqnarray*}
\end{remark}

\begin{definition}
  \label{def:hurewicz-pao}Let $X \in \mathbb{H}$, $I$ be a cosimplicial
  interval and $\sim_0$ be the relation on $\underline{\tmop{Path}} (X)$ such
  that $(n, v) \sim_0  (m, w)$ if there is a degeneracy $s$ such that $(m, w)
  = \underline{s} (n, v)$. The {\tmem{path object}} $\pao{X}$ is the
  quotient
  \[ \pao{X}  \deq  \underline{\tmop{Path}} (X) / \sim \]
  of $\underline{\tmop{Path}} (X)$ by the equivalence relation generated by
  $\sim_0$. We shall call an $\omega : \pao{X}$ {\tmem{path in $X$}} by abuse
  of language. A path $\omega : \pao{X}$ is {\tmem{constant}} if has a
  representative of degree $0$.
\end{definition}

\begin{remark}
  \label{rem:hurewicz-star}{\tmdummy}

  \begin{enumeratenumeric}
    \item The relation $\sim_0$ is reflexive and transitive yet not symmetric,
    hence $(n, v) \sim (m, w)$ iff $v$ and $w$ are connected by a zigzag of
    degeneracies. We can suppose $s$ to be elementary without restriction of
    generality.

    \item Assume $(n, v) \sim (m, w)$. If $v$ and $w$ are seen as lists (c.f.
    Remark \ref{rem:lists}), they carry the same information up to the number
    and position of occurences.

    \item We have $\pao{\astt}{} \cong \star$ since there is only the trivial
    path up to degeneracy, so the quotient collapses.
  \end{enumeratenumeric}
\end{remark}

\begin{notation}
{\tmdummy}
\begin{enumeratenumeric}
  \item Assume $(n, w) : \underline{\tmop{Path}} (X)$. We shall write
  $\reallywidetilde{(n, w)}$ for its equivalence class in $\pao{X}$.

  \item Assume $\omega : \pao{X}$. We shall write
  \begin{itemizeminus}
    \item $\omega \pic{\cdot}$ for an arbitrary representantive of $\omega$,
    in which case $\omega = \reallywidetilde{\left( \ell \left( \omega \pic{\cdot}
    \right), \omega \pic{\cdot} \right)}$;

    \item $\omega \pic{n}$ when we need to insist that $\ell \left( \omega
    \pic{\cdot} \right) = n$, in which case $\omega = \reallywidetilde{\left( n,
    \omega \pic{n} \right)}$.
  \end{itemizeminus}
\end{enumeratenumeric}
\end{notation}

\begin{lemma}
  \label{lem:zigzag}Assume $X, Y \in \mathbb{T}$ and $t :
  \underline{\tmop{Path}} (X) \rightarrow Y$. The following are equivalent
  \begin{enumerateroman}
    \item $t$ is constant on equivalence classes of $\sim$;

    \item $\left( t \circ \underline{s_i} \right) (n, w) = t (n, w)$ for all
    $n : \nno$, $(n, w) \in \underline{\tmop{Path}} (X)$ and all elementary
    degeneracies
    \begin{eqnarray*}
      s_i : \tmop{Path} (X)_n & \rightarrow & \tmop{Path} (X)_{n + 1}
    \end{eqnarray*}
  \end{enumerateroman}
\end{lemma}

\begin{proof}
  We only need to show the implication $(i \nocomma i) \Rightarrow (i)$.
  Assume $(n, v) \sim (m, w)$, so $v$ and $w$ are connected by a zigzag of
  elementary degeneracies. Reading the hypothesis $(i \nocomma i)$ from left
  to right and from right to left respectively covers the two possible cases
  encountered in a zigzag. These are the base cases for an induction on the
  length of the zigzag.
\end{proof}

\begin{remark}
  \label{rem:sing}Assume $u : \pco (X)_m$ and $v : \pco (X)_n$ such that $(m,
  u) \sim (n, v)$. \ Degeneracies being monotone, we have
  \begin{eqnarray*}
    u_{\#0} & = & v_{\#0}\\
    u_{\# (m - 1)} & = & v_{\# (n - 1)}
  \end{eqnarray*}
  hence the source and target morphisms $\src_X, \tgt_X : \pao{X}
  \rightarrow X$ given by
  \begin{eqnarray*}
    \src_X & \deq & \lambda \omega : X^{\langle I \rangle} . \omega \langle
    \cdot \rangle_{\#0}\\
    \tgt_X & \deq & \lambda \omega : X^{\langle I \rangle} . \omega \langle
    \cdot \rangle_{\# (n - 1)}
  \end{eqnarray*}
  in terms of an arbitrary representative are well-defined. Hence
  \begin{itemizeminus}
    \item $\pao{\partial^{_{} -}_X, \partial^+_X : X}{I} \rightrightarrows X$
    is an internal graph in $\mathbb{T}$;

    \item composition ``by concatenation'' $\otimes_X : \pao{X} \times \pao{X}
    \rightarrow \pao{X}$ is well-defined;

    \item the constant path morphism $\iota_X \deq \lambda x : X.
    \reallywidetilde{(0, x)} : X \rightarrow \pao{X}$ is a section of both $\src_X$
    and $\tgt_X$ so in particular $X \vartriangleleft \pao{X}$.
  \end{itemizeminus}
\end{remark}

\begin{theorem}
  {\tmdummy}

  \begin{enumeratenumeric}
    \item $\pao{X} \rightrightarrows X$ is an internal category with object
    of objects $X$, object of morphisms $\pao{X}$, composition $\otimes_X$ and
    unit $\iota_X$.

    \item There is an involution $\text{{\tmstrong{{\tmem{rev}}}}} : (-) :
    \pao{X} \rightarrow \pao{X}$ given by list reversal.

    \item The assignment \ $\pao{(-)}{I} : \mathbb{T} \rightarrow \mathbb{T}$
    is functorial, acting on morphisms by postcomposition.

    \item The morphisms $\iota_X$, $\src_X$ and $\tgt_X$ are natural in $X$.
  \end{enumeratenumeric}
\end{theorem}

\begin{proof}
  Assume $\omega : \pao{X}$.
  \begin{enumeratenumeric}
    \item By Remark \ref{rem:sing};

    \item Let $\rev \left( n, \omega \pic{n} \right)  \deq  \reallywidetilde{\left(
    n, \wp \left( \omega \pic{n} \right) \right)}$. We have
    \begin{eqnarray*}
      \left( \rev \circ \underline{s_i} \right)  \left( n, \omega \pic{n}
      \right) & = & \rev \left( n + 1, s_i \left( \omega \pic{n} \right)
      \right)\\
      & = & \reallywidetilde{\left( n + 1, \wp \left( s_i \left( \omega \pic{n}
      \right) \right) \right)}\\
      & = & \reallywidetilde{\left( n + 1, s_{n - i} \left( \wp \left( \omega
      \pic{n} \right) \right) \right)}\\
      & = & \reallywidetilde{\left( n, \wp \left( \omega \pic{n} \right) \right)}\\
      & = & \rev \left( n, \omega \pic{n} \right)
    \end{eqnarray*}
    so $\rev$ is well-defined by Lemma \ref{lem:zigzag};

    \item Let $f : X \rightarrow Y$ be a morphism in $\mathbb{T}$ and
    $\pao{f}$ be the morphism contructed by the term
    \[ \lambda \left( n, \omega \pic{n} \right) . \reallywidetilde{\left( n, \map
       \left( f, \omega \pic{n} \right) \right)} \]
    Notice that $\map \left( f, \omega \pic{n} \right)$ correspnds to $f \circ
    \omega \pic{n}$ under the isomorphism of Remark \ref{rem:lists}. Assume $0
    \leqslant i < n$. We have
    \begin{eqnarray*}
      \map \left( f, s_i \left( \omega \pic{n} \right) \right) & = & s_i
      \left( \map \left( f, \omega \pic{n} \right) \right)
    \end{eqnarray*}
    hence
    \begin{eqnarray*}
      \left( \pao{f} \circ \underline{s_i} \right) \left( n, \omega \pic{n}
      \right) & = & \pao{f} \left( n + 1, s_i \left( \omega \pic{n} \right)
      \right)\\
      & = & \reallywidetilde{\left( n + 1, \map \left( f, s_i \left( \omega \pic{n}
      \right) \right) \right)}\\
      & = & \reallywidetilde{\left( n + 1, s_i \left( \map \left( f, \omega \pic{n}
      \right) \right) \right)}\\
      & = & \reallywidetilde{\left( n, \left( \map \left( f, \omega \pic{n} \right)
      \right) \right)}\\
      & = & \pao{f} \left( n, \omega \pic{n} \right)
    \end{eqnarray*}
    so $\pao{f}$ is well-defined by Lemma \ref{lem:zigzag};

    \item obvious.
  \end{enumeratenumeric}
\end{proof}

\begin{notation}
Assume $X \in \mathbb{T}$ and $x, x' : X$. We shall write $\omega
: x \rightsquigarrow x'$ as an abbreviation for a path $\omega : \pao{X}$
such that $\src (\omega) = x$ and $\tgt (\omega) = x'$.
\end{notation}

\begin{definition}
  A face filtration is a sequence $d = (d^{(i)})_{0 < i \leqslant n}$ of
  morphisms where $d^{(n)} = \tmop{id}$ and $d = (d^{(i)})_{0 < i < n}$ are
  faces in
  \[ \pco (X)_{\mu_d (0)} \longleftarrowlim^{d^{(1)}} \pco (X)_{\mu_d (1)}
     \longleftarrowlim^{d^{(2)}} \cdots \longleftarrowlim^{d^{(n - 2)}} \pco
     (X)_{\mu_d (n - 1)} \longleftarrowlim^{d^{(n - 1)}} \pco (X)_{\mu_d (n)}
  \]
  with $\mu_d : \{ 0, \cdots, n \} \rightarrow \nno$.
\end{definition}

\begin{notation}
  $d \flatt{}{i}  \deq d^{(i)} \circ d^{(i + 1)} \circ \cdots \circ
  d^{(n)}$ for a face filtration $(d^{(i)})_{0 < i \leqslant n}$.
\end{notation}

\begin{remark}
  We have in particular
  \begin{eqnarray*}
    \flatt{d}{n} (w) & = & d^{(n)} (w)\\
    & = & w\\
    \flatt{d}{1} (w) & = & (d^{(1)} \circ \cdots \circ d^{(n - 1)}) (w)
  \end{eqnarray*}
\end{remark}

\begin{example}
  {\tmdummy}

  \begin{enumeratenumeric}
    \item Assume $n : \nno$. The {\tmem{tail filtration}} $\tmop{tail}_n =
    (\tmop{tail}_n^{(i)})_{0 \leqslant i < n}$ is given \ by the {\tmem{first
    face}}
    \begin{eqnarray*}
      \tmop{tail}_n^{(i)} & \deq & d_0
    \end{eqnarray*}
    while the {\tmem{prefix filtration}} $\tmop{pref}_n =
    (\tmop{pref}_n^{(i)})_{0 \leqslant i \leqslant n}$ is given by the
    {\tmem{last face}}
    \begin{eqnarray*}
      \tmop{pref}_n^{(i)} & \deq & d_i
    \end{eqnarray*}
    at $0 < i < n$ respectively.

    \item Assume $m, n : \nno$. The {\tmem{chop filtration}} $\tmop{chop}_{m,
    n} = (\tmop{chop}^{(i)}_{m, n})_{0 \leqslant i < n}$ is given by
    \begin{eqnarray*}
      \mu_{\tmop{chop}_{m, n}} (n) & \deq & 2 n + m\\
      \tmop{chop}_{m, n}^{(i)} & \deq & \tmop{tail}_n^{(i)} \circ
      \tmop{pref}_n^{(i - 1)}\\
      & = & \tmop{pref}_n^{(i)} \circ \tmop{tail}_n^{(i - 1)}
    \end{eqnarray*}
    We have in particular $\mu_{\tmop{chop}_{m, n}} (0) = m$.
  \end{enumeratenumeric}
  We have $\flatt{\tmop{pref}}{1} (w) = [w_{\#0}]$ a constant path and
  similarly for the tail filtration as well as for the chop filtration.
\end{example}

\begin{definition}
  Assume a cosimplicial interval $I$.
  \begin{enumeratenumeric}
    \item assume a face filtration $d = (d^{(i)})_{0 < i \leqslant n}$ and $w
    : \pco (X)_{\mu_d (n)}$. We shall call the list
    \[ \tmop{eval} (d, w)  \deq \left[ \reallywidetilde{\left( \flatt{\mu_d (0),
       d}{1} (w) \right)} ; \cdots ; \reallywidetilde{\left( \mu_d (n - 1),
       \flatt{d}{n} (w) \right)} \right] : \tmop{List}_n \left( \pao{X}
       \right) \]
    evaluation of $d$ at $w$;

    \item The interval $I$ has the {\tmem{Hurewicz property}} if for any face
    filtration $d = (d^{(i)})_{0 < i \leqslant n}$ and any $w : \pco
    (X)_{\mu_d (n)}$ we have
    \begin{eqnarray*}
      \tmop{eval} (d, w) & \in & {\pao{X}}^{I_n}
    \end{eqnarray*}
    That is, any evaluation of a face filtration is a (rigid) path of paths.
  \end{enumeratenumeric}
\end{definition}

\begin{theorem}
  $\Delta \tmmathbf{2} \in \eff$ is a Hurewicz interval.
\end{theorem}

\begin{proof}
  $\Delta \tmmathbf{2}$ is cosimplicial (c.f. Remark \ref{ex:cosimp}). To see
  that it verifies the Hurewicz property, recall that the set underlying an
  exponential $(Y, \approx)^{(X, \approx)}$ in {\eff} is $\mathcal{P}
  (\mathbb{N})^{X \times Y}$, while existence is the set of (encoded) 4-tuples
  $\langle a, b, c, d \rangle$ asserting that an element of $\mathcal{P}
  (\mathbb{N})^{X \times Y}$ is a functional relation. Given a face $d :
  X^{I_m} \rightarrow X^{I_n}$ we have
  \begin{eqnarray*}
    E (d (w)) & \subseteq & E (w) \quad (\ast)
  \end{eqnarray*}
  since precomposing with a coface only retains some of the original 4-tuples.
  On the other hand, given a degeneracy $s : X^{I_n} \rightarrow X^{I_m}$ we
  have
  \begin{eqnarray*}
    E (s (w)) & = & E (w) \quad (\ast \ast)
  \end{eqnarray*}
  since precomposing with a codegeneracy does not add any new 4-tuple. Assume
  \begin{itemizeminus}
    \item a face filtration $d = (d^{(i)})_{0 < i \leqslant n}$;

    \item $w : \pco (X)_{\mu_d (n)}$;

    \item an ordered subset $\{ s_0, \cdots, s_r \} \subset \{ 0, \cdots, n
    \}$.
  \end{itemizeminus}
  Let $w_k  \deq  \flatt{d}{s_k} (w)$ and $n_k  \deq \mu_d (s_k) $for $0
  \leqslant k \leqslant r$. Applying $(\ast)$ yields the ascending chain of
  inclusions of realising sets
  \[ E (w_0) \subset E (w_1) \subset \cdots \subset E (w_r) \quad (\ast \ast
     \ast) \]
  but only if we consider those paths {\tmem{in isolation}}, that is if we
  consider
  \begin{eqnarray*}
    & w_k : & \tmop{Path} (X)_{n_k}
  \end{eqnarray*}
  However, these inclusions do not hold anymore if we inject the $w_k$'s into
  $\underline{\tmop{Path}} (X)$, that is if we consider
  \begin{eqnarray*}
    (n_k, w_k) & : & \underline{\tmop{Path}} (X)
  \end{eqnarray*}
  For existence becomes then a set of 5-tuples $\langle n_k, a, b, c, d
  \rangle$ with the degree tag $n_k$ as additional datum. Assume
  \begin{eqnarray*}
    \langle a, b, c, d \rangle & \in & E (w_0)
  \end{eqnarray*}
  By $(\ast \ast \ast)$ we have
  \begin{eqnarray*}
    \langle n_k, a, b, c, d \rangle & \in & E (n_k, w_k)
  \end{eqnarray*}
  for all $k \in \{ 0, \cdots, r \}$. But for any $k \in \{ 0, \cdots, r \}$
  we can construct a degeneracy $s_{k, r} : I^{n_k} \rightarrow I^{n_r}$ (for
  instance a ``one-sided stuffing'') such that the endomorphism
  $\underline{s_{k, r}} : \underline{\tmop{Path}} (X) \rightarrow
  \underline{\tmop{Path}} (X)$ adjusts the degree tag to $n_r$. By $(\ast
  \ast)$ $\underline{s_{k, r}}$ does not alter the other coordinates of a
  realiser, so we have
  \begin{eqnarray*}
    \langle n_r, a, b, c, d \rangle & \in & E \left( \underline{s_{k, r}}
    (n_k, w_k) \right)
  \end{eqnarray*}
  This entails
  \begin{eqnarray*}
    \langle n_r, a, b, c, d \rangle & \in & \bigcup_{(\deg (w), w) \sim (n_k,
    w_k)} E (\deg (w), w)
  \end{eqnarray*}
  for all $k \in \{ 0, \cdots, r \}$. Hence
  $\tmop{eval} (d, w) \in {\pao{X}}^{I_n}$ by Remark \ref{rem:rpath}.
\end{proof}

Assume a topos $\mathbb{H}$ with NNO equipped with a Hurewicz interval $I$.

\begin{definition}
  \label{def:htopy}Let $f, g : X \rightarrow Y$ be morphisms. A
  {\tmem{homotopy}} $H : f \htpy g$ from $f$ to $g$ is given by a commuting
  diagram

  \begin{center}
    {\includegraphics[trim = {38 0 0 30}]{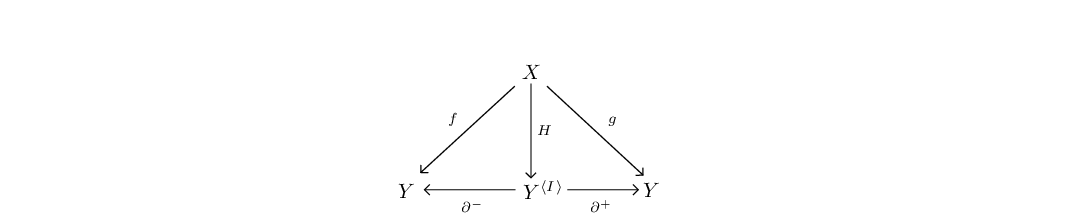}}
  \end{center}

  {\noindent}$H$ is {\tmem{constant}} on a subobject $X' \vartriangleleft X$
  provided $H (x) = \iota (x)$ for any $x : X$ such that $x \in X'$.
\end{definition}

\begin{remark}
  \label{rem:homot-path}A homotopy $H : f \looparrowright g$ informs us that
  for any $x : X$ there is a path $H (x) : \pao{Y}$ such that $f (x) = \src
  (\omega_x)$ and $g (x) = \tgt (\omega_x)$.
\end{remark}

\begin{definition}
  A {\tmem{homotopy equivalence}} is a morphism $u : X \rightarrow Y$ which
  has homotopy inverse.
\end{definition}

\begin{definition}
  Assume a category $\mathbb{C}$. A class of morphisms $\mathcal{A} \subset
  \mathbb{C}_1$ has
  \begin{enumeratenumeric}
    \item the {\tmem{3-for-2 property}} if for any factorisation $a_3 = a_2
    \circ a_1$ the membership $a_i, a_j \in \mathcal{A}$ for some $\{ i, j \}
    \subset \{ 1, 2, 3 \}$ entails $a_k \in \mathcal{A}$ for $k \in \{ 1, 2, 3
    \} \backslash \{ i, j \}$;

    \item the {\tmem{6-for-2 property}} if given morphisms $A \rightarrowlim^u
    B \rightarrowlim^v C \xrightarrow{w} D$, $v \circ u, w \circ v \in
    \mathcal{A}$ entails
    \[ u, v, w, w \circ v \in \mathcal{A} \]
  \end{enumeratenumeric}
\end{definition}

\begin{remark}
  \label{rem:three-for-two}Homotopy equivalences structurally verify 3-for-2
  and {\tmem{weak invertibility}}, that is given morphisms $A \rightarrowlim^u
  B \rightarrowlim^v C \xrightarrow{w} D$ with $v \circ u$ and $w \circ v$
  homotopy equivalences the morphism $v$ has to be a homotopy equivalence.
  This entails that they verify 6-for-2 {\cite{dwyer2005homotopy}}.
\end{remark}

\begin{remark}
  \label{rem:vac}Assume $\omega : \pao{X}$ represented by $\omega \pic{n}$ and
  $0 \leqslant i < n$. Given the elementary degeneracy $s_i : X^{I_n}
  \rightarrow X^{I_{n + 1}}$ we have
  \begin{eqnarray*}
    \tmop{eval} \left( \tmop{pref}_n, \omega \pic{n} \right)_{\#j} & = &
    \reallywidetilde{\left(n-j,\flatt{\tmop{pref}}{j} \left( \omega \pic{n} \right) \right)}\\
    \tmop{eval} \left( \tmop{pref}_{n + 1}, s_i \left( \omega \pic{n} \right)
    \right)_{\#j} & = & \left\{\begin{array}{lll}
      \reallywidetilde{\left(n-j,\flatt{\tmop{pref}}{j} \left( \omega \pic{n} \right) \right)} &  & 0
      \leqslant j \leqslant i\\
      \reallywidetilde{\left(n-j+1,\flatt{\tmop{pref}}{j - 1} \left( \omega \pic{n} \right) \right)} &
      & i + 1 \leqslant j < n + 1
    \end{array}\right.
  \end{eqnarray*}
  hence
  \begin{eqnarray*}
    \tmop{eval} \left( \tmop{pref}_n, \omega \pic{n} \right) & \sim &
    \tmop{eval} \left( \tmop{pref}_{n + 1}, s_i \left( \omega \pic{n} \right)
    \right)
  \end{eqnarray*}
  so the morphism $C_n : \pao{X} \rightarrow {\pao{X}}^{\langle I \rangle}$ constructed by the
  term $\lambda \omega : \pao{X} . \reallywidetilde{\tmop{eval} \left(
  \tmop{pref}_n, \omega \pic{\cdot} \right)}$ is well-defined by Lemma
  \ref{lem:zigzag}. It follows that everything in sight commutes in the
  diagram

  \begin{center}
    \includegraphics[trim = {38 0 0 40}]{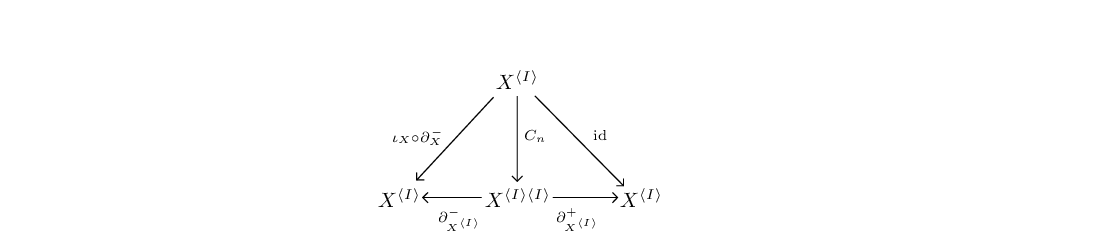}
  \end{center}

  {\noindent}The homotopy $C_n : \iota_X \circ \src_X \looparrowright
  \tmop{id}_{\pao{X}}$ is called {\tmem{contracting homotopy}}.
\end{remark}

\section{Hurewicz fibrations}\label{sec:hurewicz-fibs}

\begin{definition}
  A section $h$ of the canonical morphism $\left\langle \pao{p}, \src_E
  \right\rangle$ in

  \begin{center}
    \includegraphics[trim = {38 0 0 30}]{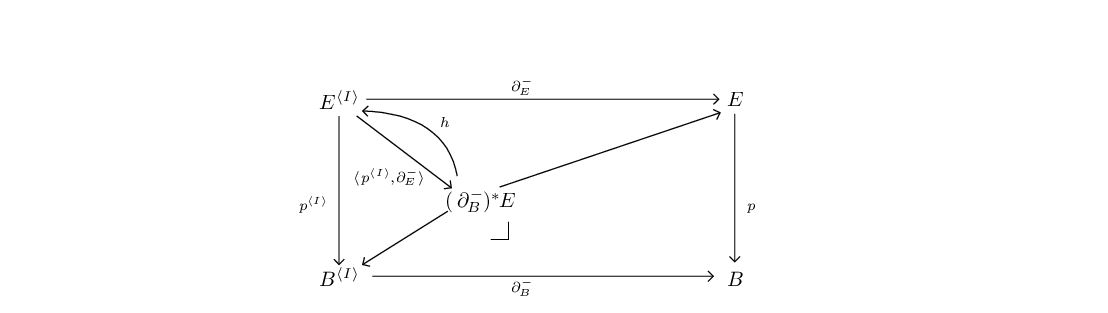}
  \end{center}

  {\noindent}is called {\tmem{connection}} (for $p$) if it preserves constant
  paths. A morphism which admits a connection is called {\tmem{(Hurewicz)
  fibration}}.
\end{definition}

\begin{remark}
  A fibration $p : E \rightarrow B$ is thus a morphism with a strong path
  lifiting property: for any path $\omega : b \rightsquigarrow b'$ in $B$ and
  any $e : E$ such that $p (e) = b$ there is a path $\varpi$ in $X$ such that
  $\pao{p} (\varpi) = \omega$ along with an explicit construction of one such
  lift. We call an $(\omega, e) : \left( \src_B \right)^{\ast} E$
  {\tmem{initial datum}} for $p$.
\end{remark}

\begin{proposition}
  \label{prop:hurewicz-compo-stable}Fibrations are stable under composition.
\end{proposition}

\begin{remark}
  The operation of pairing or ``zipping'' paths $\pao{X} \times \pao{Y}
  \rightarrow \pao{(X \times Y)}$ is in general not well-defined. However it
  is obviously the case when one of the arguments is constant, so we have the
  morphism
  \begin{eqnarray*}
    \pad_l : X \times \pao{Y} & \rightarrow & \pao{(X \times Y)}
  \end{eqnarray*}
  constructed by the term $\lambda (x, \omega) : X \times \pao{Y} .
  \reallywidetilde{\map (\lambda y : Y. (x, y)) \left( \omega \pic{\cdot}
  \right)}$ and the morphism
  \begin{eqnarray*}
    \pad_r : \pao{X} \times Y & \rightarrow & \pao{(X \times Y)}
  \end{eqnarray*}
  constructed by the term $\lambda (\omega, y) : \pao{X} \times Y.
  \reallywidetilde{\map (\lambda x : X. (x, y)) \left( \omega \pic{\cdot}
  \right)}$.
\end{remark}

\begin{proposition}
  Projections from products are fibrations.
\end{proposition}

\begin{proof}
  Assume $(x, y) \in X \times Y$ and $\omega : \pao{X}$. There is the obvious
  path $\bar{\omega} : (x, y) \rightsquigarrow (x', y)$ constant in the second
  coordinate. The term
  \[ \lambda (\omega, (x, y)) : \left( \src_X \right)^{\ast} (X \times Y) .
     \text{{\tmstrong{pad}}}_r (\omega, y) \]
  constructs a connection for $\pi_0$. Similarly for $\pi_1$.
\end{proof}

\begin{proposition}
  \label{prop:src-tgt}$\left\langle \src, \tgt \right\rangle : \pao{X}
  \rightarrow X \times X$ is a fibration for any $X \in \mathbb{H}$.
\end{proposition}

\begin{proof}
  Assume $\omega : \pao{X}$ and $\kappa : \pao{(X \times X)}$ $(\omega,
  \kappa)$ is an initial datum for $\piq{\partial_X^-}{\partial_X^+}$, that is
  such that $\piq{\partial_X^-}{\partial_X^+} (\omega) = \partial_{X \times
  X}^- (\kappa)$. We have
  \begin{eqnarray*}
    \kappa_0 \pic{n} & \deq & \map (\pi_0) \left( \kappa \pic{n} \right) \in
    X^{I_n}\\
    \kappa_1 \pic{n} & \deq & \map (\pi_1) \left( \kappa \pic{n} \right) \in
    X^{I_n}\\
    \tmop{unwind} (\kappa, \omega) & \deq & \rev \left( \kappa_0 \pic{n}
    \right) \otimes \omega \pic{m} \otimes \kappa_1 \pic{n} \in X^{I_{2
    \nocomma n + m}}
  \end{eqnarray*}
  hence
  \[ \tmop{eval} \left( \tmop{chop}_{m, n}, \tmop{unwind} \left( \omega
     \pic{m}, \kappa \pic{n} \right) \right) \in {\pao{X}}^{I_n} \]
  by the Hurewicz property.
  We need to show that this term is constant on equivalence classes. Notice that by construction
  \begin{eqnarray*}
    \tmop{eval} \left( \tmop{chop}_{m, n}, \tmop{unwind} \left( \omega
    \pic{m}, \kappa \pic{n} \right) \right)_{\#0} & = & \omega
  \end{eqnarray*}
  for any representant $\omega \pic{m}$ of $\omega$. Assume $0 < j \leqslant
  n$ and let
  \begin{eqnarray*}
    l_n^{(j)} & \deq & \reallywidetilde{\left(n-j,\rev \left( (\tmop{pref}_n)_{\flat}^{(j)} \left(
    \kappa_0 \pic{n} \right) \right) \right)}\\
    r_n^{(j)} & \deq & \reallywidetilde{\left(n-j,(\tmop{pref}_n)_{\flat}^{(j)} \left( \kappa_1 \pic{n}
    \right)\right)}
  \end{eqnarray*}
  We then have
  \begin{eqnarray*}
    \tmop{eval} \left( \tmop{chop}_{m, n}, \tmop{unwind} \left( \omega
    \pic{m}, \kappa \pic{n} \right) \right)_{\#j} & = & l_n^{(j)}
    \otimes \omega \otimes r_n^{(j)}
  \end{eqnarray*}
  Assume $0 \leqslant i < n$ and the elementary degeneracy $s_i : X^{I_n}
  \rightarrow X^{I_{n + 1}}$. Let
  \begin{eqnarray*}
    l^{(i, j)}_{n + 1} & \deq & \left\{\begin{array}{lll}
      l^{(j)}_n &  & 0 < j \leqslant i\\
      \reallywidetilde{\left( n-j+1, s_{j - i - 1} \left( \rev \left( (\tmop{pref}_n)_{\flat}^{(j)} \left(
      \kappa_0 \pic{n} \right) \right) \right) \right)} &  & i < j \leqslant n + 1
    \end{array}\right.\\
    r^{(i, j)}_{n + 1} & \deq & \left\{\begin{array}{lll}
      r^{(j)}_n &  & 0 < j \leqslant i\\
      \reallywidetilde{\left( n-j+1, s_i \left( (\tmop{pref}_n)_{\flat}^{(j)} \left( \kappa_1 \pic{n} \right)
      \right)  \right)} &  & i < j \leqslant n + 1
    \end{array}\right.
  \end{eqnarray*}
  We have
  \begin{eqnarray*}
    \tmop{eval} \left( \tmop{chop}_{m, n + 1}, \tmop{unwind} \left( \omega
    \pic{m}, s_i \left( \kappa \pic{n} \right) \right) \right)_{\#j} & = &
    l_{n + 1}^{(i, j)} \otimes \omega \otimes r_{n +
    1}^{(i, j)}\\
    & = & \left\{\begin{array}{lll}
      l_n^{(j)} \otimes \omega \otimes r_n^{(j)} &  &
      0 < j \leqslant i\\
      l_n^{(j - 1)} \otimes \omega \otimes r_n^{(j -
      1)} &  & i + 1 \leqslant j \leqslant n
    \end{array}\right.
  \end{eqnarray*}
  so
  \begin{eqnarray}
    \tmop{eval} \left( \tmop{chop}_{m, n}, \tmop{unwind} \left( \omega
    \pic{m}, \kappa \pic{n} \right) \right) & \sim & \tmop{eval} \left(
    \tmop{chop}_{m, n + 1}, \tmop{unwind} \left( \omega \pic{m}, s_i \left(
    \kappa \pic{n} \right) \right) \right) \nonumber
  \end{eqnarray}
  hence the term
  \[ \lambda (\kappa, \omega) : (\partial^-_{X \times X})^{\ast} \pao{X} .
     \reallywidetilde{\tmop{eval} \left( \tmop{chop}_{m, n}, \tmop{unwind} \left(
     \omega \pic{\cdot}, \kappa \pic{\cdot} \right) \right)} \]
  is well-defined by Lemma \ref{lem:zigzag}. The section
  $h_{\piq{\src_X}{\tgt_X}} : (\partial^-_{X \times X})^{\ast} \pao{X}
  \rightarrow {\pao{X}}^{\langle I \rangle}$ constructed by this term preserves constant
  paths.
\end{proof}

\begin{corollary}
  The source map $\src_X : \pao{X} \rightarrow X$ and the target map
  $\tgt_X : \pao{X} \rightarrow X$ are fibrations for any $X \in
  \mathbb{H}$.
\end{corollary}

\section{Weak Factorisation System}\label{sec:wfs}

\begin{definition}
  $X$ is a {\tmem{strong deformation retract}} of $Y$ if there is a morphism
  $e : X \rightarrow Y$ admitting a retraction $r : Y \rightarrow X$ such that
  there is a homotopy $H : e \circ r \looparrowright \tmop{id}_Y$ constant on
  $X$. We call the split epi $r$ {\tmem{strong deformation retraction }}and
  the split mono $e$ {\tmem{strong deformation insertion}}, respectively.
\end{definition}

\begin{notation}
  We shall write $\tmop{SDI}$ for the class of strong deformation
  insertions and $\mathcal{H}$ for the class of fibrations
\end{notation}

\begin{remark}
  A strong deformation insertion is a homotopy equivalence.
\end{remark}

\begin{definition}
  Let $f : X \rightarrow Y$ be a morphism in $\mathbb{H}$. The object $M_f$
  given by the pullback

  \begin{center}
    \includegraphics[trim = {38 0 0 30}]{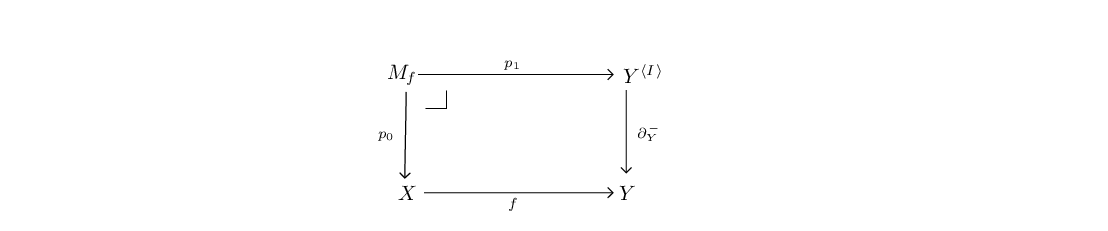}
  \end{center}

  {\noindent}is called $f$'s {\tmem{mapping track}}.
\end{definition}

\begin{remark}
  $M_f = \left\{ (x, \omega) : X \times \pao{Y} |f (x) = \src_Y (\omega)
  \right\}$ is the object of paths that begin in the image of $f$.
\end{remark}

\begin{theorem}
  \label{th:facto}A morphism $f : X \rightarrow Y$ factors through the mapping
  track as a strong deformation insertion followed by a fibration.
\end{theorem}

\begin{proof}
  Assume

  \begin{center}
    \includegraphics[trim = {-10 20 0 60}]{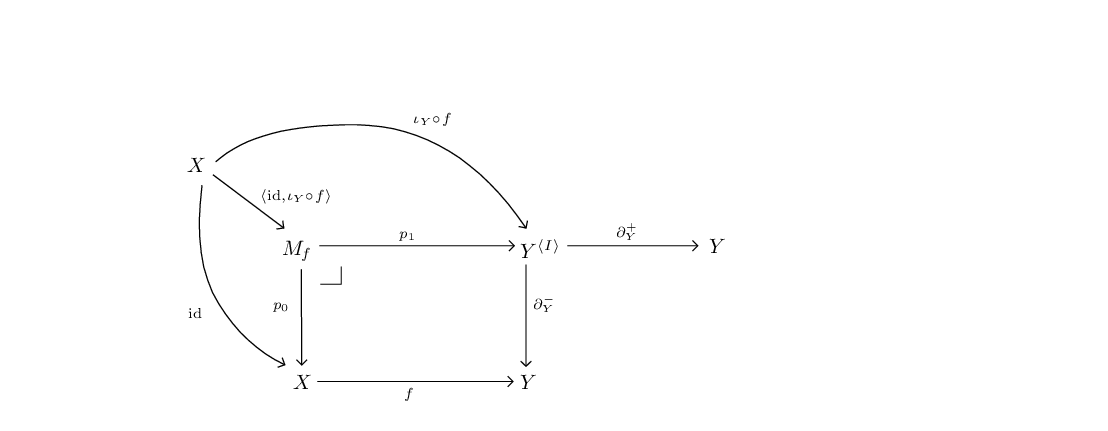}
  \end{center}

  {\noindent}We have
  \begin{eqnarray*}
    f \circ \tmop{id} & = & \src_Y \circ \iota_Y \circ f
  \end{eqnarray*}
  since $\iota_Y$ is a section of $\src_Y$, while
  \begin{eqnarray*}
    \left( \tgt_Y \circ p_1 \right) \circ \langle \tmop{id}_X, \iota_Y \circ f
    \rangle & = & \tgt_Y \circ \iota_Y \circ f\\
    & = & f
  \end{eqnarray*}
  since $\iota_Y$ is a section of $\tgt_Y$. This is the factorisation we seek
  since
  \begin{enumeratenumeric}
    \item $p_0$ is a retraction of $\langle \tmop{id}_X, \iota_Y \circ f
    \rangle$ by construction. The term
    \[ \lambda (x, \omega) : M_f . \pad_l (x, C_n (\omega)) \]
    constructs a contracting homotopy $H_n : \langle \tmop{id}_X, \iota_Y
    \circ f \rangle \circ p_0 \looparrowright \tmop{id}$ which is constant on
    $X$ by construction (c.f. Remark \ref{rem:vac});

    \item assume $(x, \omega) : X \times_Y \pao{Y}$ and $\kappa : \pao{Y}$
    such that $\tgt (\omega) = \src \kappa$. The term
    \[ \lambda ((x, \omega), \kappa) : M_f . \pad_l \left( x,
       \reallywidetilde{\tmop{eval} \left( \tmop{pref}_n, \omega \pic{\cdot} \otimes
       \kappa \pic{\cdot} \right)} \right) \]
    is well-defined and constructs a connection for $\tgt_Y \circ p_1$. The
    argument here is essentially the simpler half of the one used in the proof
    of Proposition \ref{prop:src-tgt}.
  \end{enumeratenumeric}

\end{proof}

\begin{proposition}
  \label{prop:lift1}$\tmop{SDI} \subset \llp{\mathcal{H}}$ and
  $\rlp{\tmop{SDI}} \supset \mathcal{H}$.
\end{proposition}

\begin{proof}
  We claim that there is a lift $d : A \rightarrow E$ in any commuting diagram

  \begin{center}
    \includegraphics[trim = {70 10 0 35}]{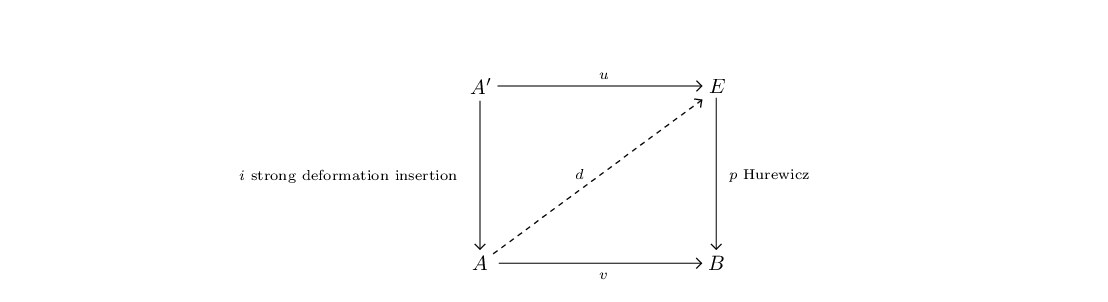}
  \end{center}

  {\noindent}Assume $r : A \rightarrow A'$ a strong deformation retraction
  with section $i$. Assume $a : A$. We have by hypothesis a homotopy $H : i
  \circ r \looparrowright \tmop{id}_A$, so
  \begin{eqnarray*}
    (i \circ r) (a) & = & \src_A (H (a))\\
    a & = & \tgt_A (H (a))
  \end{eqnarray*}
  hence
  \begin{eqnarray*}
    (v \circ i \circ r) (a) & = & \src_B \left( \pao{v} (H (a)) \right) \quad
    (\ast)\\
    v (a) & = & \tgt_B \left( \pao{v} (H (a)) \right) \quad (\ast \ast)
  \end{eqnarray*}
  by naturality. We also have
  \begin{eqnarray*}
    (v \circ i \circ r) (a) & = & (p \circ u \circ r) (a) : B \quad (\ast \ast
    \ast)
  \end{eqnarray*}
  by hypothesis, so in particular
  \begin{eqnarray*}
    p ((u \circ r) (a)) & = & \src_B \left( \pao{v} (H (a)) \right) \quad
    (\ast \ast \ast)  \text{and } (\ast)
  \end{eqnarray*}
  so $(v (H (a)), (u \circ r) (a))$ is an initial datum for $p$. Assume a
  connection $h_p : \left( \src_B \right)^{\ast} E \rightarrow \pao{E}$ for
  $p$ and let
  \begin{eqnarray*}
    d & \deq & \lambda a : A. \tgt_E \left( h_p \left( \pao{v} (H (a)), (u
    \circ r) (a) \right) \right)
  \end{eqnarray*}
  We have
  \begin{eqnarray*}
    (p \circ d) (a) & = & p (d (a))\\
    & = & \tgt_B \left( \pao{v} (H (a)) \right) \quad d (a) \text{ lift of
    $\pao{v} (H (a))$ and naturality}\\
    & = & v (a) \hspace{6em} (\ast \ast)
  \end{eqnarray*}
  Assume now $a : A$ such that $a \in A'$. We have
  \begin{eqnarray*}
    H (i (a)) & = & \iota_A (i (a)) \hspace{5em} \text{$H$ is constant on $A'$
    by hypothesis}\\
    \Rightarrow \text{\quad} \pao{v} (H (i (a))) & = & \iota_B (i (a))\\
    \Rightarrow \text{\quad$h_p \left( \pao{v} (H (i (a))), u (a) \right)$} &
    = & h_p (\iota_B (i (a)), u (a))\\
    & = & \iota_E (u (a)) \quad \text{{\hspace{4em}}connections preserve
    constant paths}
  \end{eqnarray*}
  Hence
  \begin{eqnarray*}
    (d \circ i) (a) & = & d (i (a))\\
    & = & \tgt_E \left( h_p \left( \pao{v} (H (i (a))), (u \circ r) (i (a))
    \right) \right)\\
    & = & \tgt_E \left( h_p \left( \pao{v} (H (i (a))), u (a) \right)
    \right)\\
    & = & u (a)
  \end{eqnarray*}
\end{proof}

\begin{proposition}
  \label{prop:fibrant}Any object $X \in \mathbb{H}$ is fibrant.
\end{proposition}

\begin{proof}
  A pullback over $\astt$ is just a product while $\pao{\astt}$ is terminal,
  so we get the diagram

\begin{center}
  \includegraphics[trim = {38 -5 0 30}]{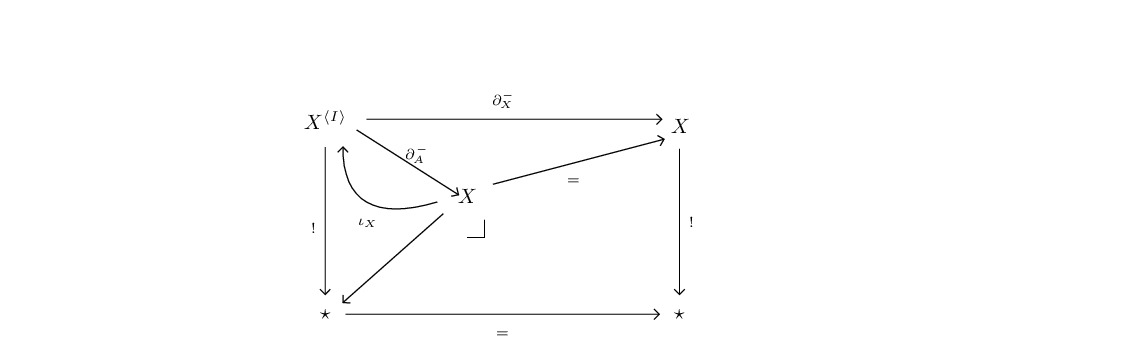}
\end{center}

  {\noindent}The outer square is the naturality square for $\src$ at $!_X$,
  the inner square is the inscribed pullback while $\src_X$ is the canonical
  morphism to the pullback. This morphism admits the section $\iota_X$.
\end{proof}

\begin{proposition}
  \label{prop:lift2}$\tmop{SDI} \supset \llp{\mathcal{H}}$.
\end{proposition}

\begin{proof}
  Assume $j : A \rightarrow X$ such that $j \in \llp{\mathcal{H}}$. Since
  every object is fibrant (c.f. Proposition \ref{prop:fibrant}), everything in
  sight commutes in the following diagram

  \begin{center}
    \includegraphics[trim = {38 15 0 40}]{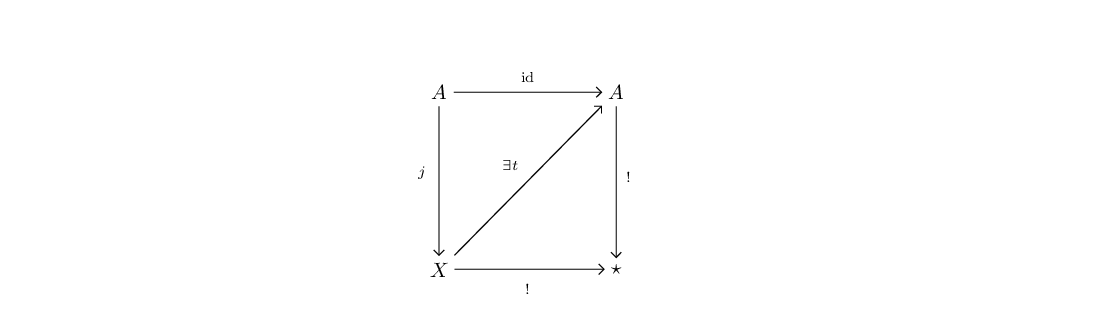}
  \end{center}

  {\noindent}so $j$ admits a retraction. Since $\piq{\src}{\tgt}$ is Hurewicz,
  \ everything in sight commutes in the following diagram

  \begin{center}
    \includegraphics[trim = {38 15 0 40}]{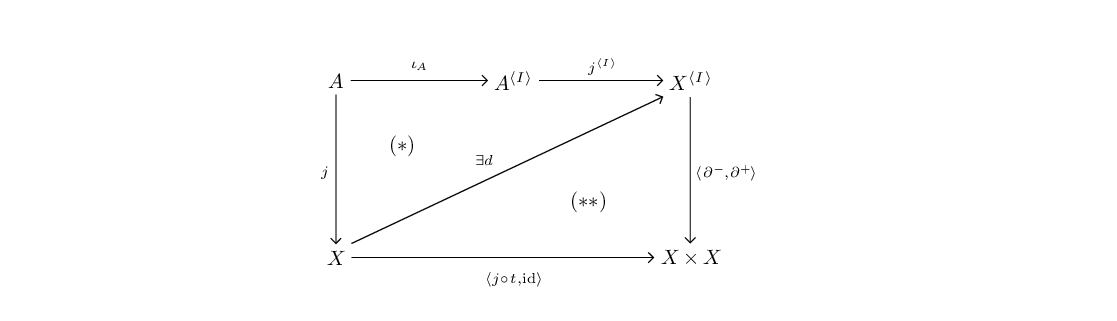}
  \end{center}

  {\noindent}so by $(\ast \ast)$ there is a homotopy $H : j \circ t
  \looparrowright \tmop{id}_X$ given by $H (x)  \deq d (x)$. This homotopy is
  constant on $A$ by $(\ast)$.
\end{proof}

\begin{proposition}
  \label{prop:lift3}$\rlp{\tmop{SDI}} \subset \mathcal{H}$.
\end{proposition}

\begin{proof}
  Assume $p : E \rightarrow B$ such that $p \in \rlp{\tmop{SDI}}$. We have the
  factorisation $p = f \circ c$ with $c \in \tmop{SDI}$ and $f \in
  \mathcal{H}$ (c.f Theorem \ref{th:facto}), hence everything in sight
  commutes in the following diagram

  \begin{center}
    \includegraphics[trim = {38 5 0 40}]{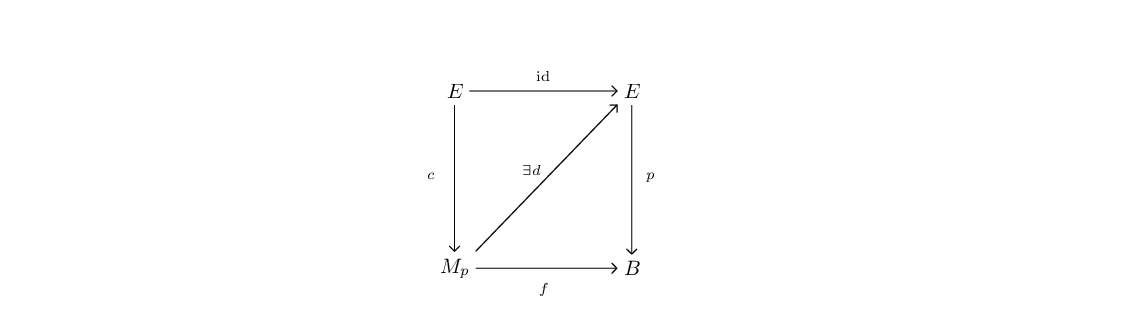}
  \end{center}

  {\noindent}Assume $(\omega, e) : \left( \src_B \right)^{\ast} E$. Assume a
  connection $h_f : \left( \src_B \right)^{\ast} M_p \rightarrow \pao{M_p}$
  for $f$. We have
  \begin{eqnarray*}
    f (c (e)) & = & p (d (c (e)))\\
    & = & p (e)
  \end{eqnarray*}
  so $(\omega, c (e)) : \left( \src_B \right)^{\ast} M_p$ is an initial datum
  for $h_f$, hence the term
  \[ \lambda (e, \omega) : \left( \src_B \right)^{\ast} E.d (h_f (c (e),
     \omega)) \]
  constructs a connection for $p$.
\end{proof}

\begin{corollary}
  \label{cor:pullback}Fibrations are closed under pullbacks and retracts.
\end{corollary}

\begin{definition}
  Assume a category $\mathbb{C}$. Classes of morphisms $\mathcal{L},
  \mathcal{R} \subset \mathbb{C}_1$ form a {\tmem{weak factorisation system}}
  provided
  \begin{enumerateroman}
    \item every morhism $f \in \mathbb{C}_1$ factors as $f = r \circ l$ with
    $l \in \mathcal{L}$ and $r \in \mathcal{R}$;

    \item $\mathcal{L}= \llp{\mathcal{R}}$;

    \item $\rlp{\mathcal{L}} =\mathcal{R}$.
  \end{enumerateroman}
\end{definition}

\begin{theorem}
  The classes $\tmop{SDI}$ and $\mathcal{H}$ form a weak factorisation system.
\end{theorem}

\begin{proof}
  The factorisation is given by Theorem \ref{th:facto} while the lifting
  conditions are a consequence of Propositions \ref{prop:lift1},
  \ref{prop:lift2} and \ref{prop:lift3}.
\end{proof}

\section{Category of Fibrant Objects}\label{sec:cfo}

\begin{proposition}
  \label{prop:triv-pullback}Trivial fibrations are closed under pullback.
\end{proposition}

\begin{proof}
  Assume $p : E \rightarrow B$ a fibration witnessed by connection $h_p$.
  Assume $f : A \rightarrow B$. Now $p_0  \deq f^{\ast} p$ is a fibration
  (c.f. Corollary \ref{cor:pullback}), so we only need to establish that it is
  a homotopy equivalence. Assume a homotopy inverse $u : B \rightarrow E$ of
  $p$ witnessed by homotopies
  \begin{eqnarray*}
    H : p \circ u & \looparrowright & \tmop{id}_B\\
    K : u \circ p & \looparrowright & \tmop{id}_E
  \end{eqnarray*}
  Assume $a : A$. We have

  \begin{center}
    \includegraphics[trim = {38 0 0 40}]{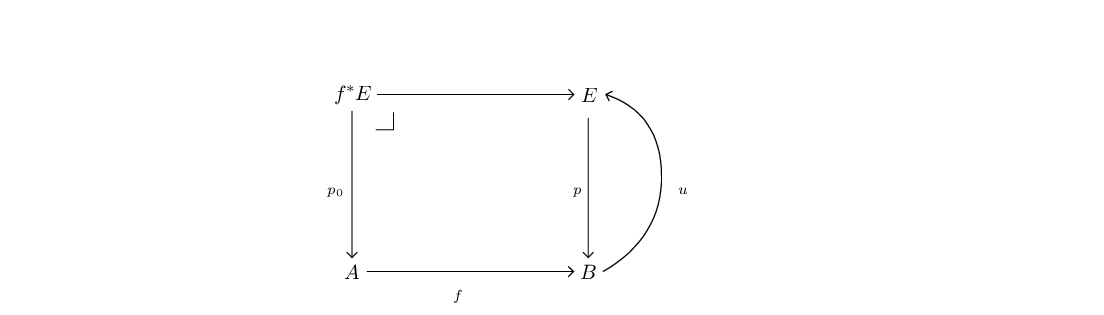}
  \end{center}

  {\noindent}But $p$ is Hurewicz, so we have the lift $h_p (u (f (a)), H (f
  (a)))$. Let
  \[ e_a  \deq  \tgt (h_p (u (f (a)), H (f (a)))) \]
  As $p (e_a) = f (a)$, the term $\lambda a : A. (a, e_a)$ constructs a
  section $u_0 : A \rightarrow f^{\ast} E$ of $p_0$. We claim that $u_0$ is a
  deformation insertion. Assume $(a, e) : f^{\ast} E$. We have $u_0 (p_0 (a,
  e)) = u_0 (a) = (a, e_a)$ and paths
  \begin{eqnarray*}
    K (e) : u (p (e)) & \rightsquigarrow & e\\
    K (e_a) : u (p (e_a)) & \rightsquigarrow & e_a
  \end{eqnarray*}
  But $p (e) = f (a)$ by hypothesis and $p (e_a) = f (a)$ by construction so
  $u (p (e)) = u (p (e_a))$. Hence the term
  \[ \lambda (a, e) : f^{\ast} E. \text{{\tmstrong{pad}}}_l \left( \iota_A
     (a), \rev (K (e_a)) \otimes K (e) \right) \]
  constructs a homotopy $H' : u_0 \circ p_0 \looparrowright
  \tmop{id}_{f^{\ast} E}$.
\end{proof}

\begin{proposition}
  \label{prop:constant-path-sd}The constant path morphism $\iota_X : X
  \rightarrow \pao{X}$ is a strong deformation insertion.
\end{proposition}

\begin{proof}
  By remark \ref{rem:vac}.
\end{proof}

\begin{remark}
  \label{rem:diagonal-facto}The diagonal factors through $\pao{X}$ as a
  homotopy equivalence (a strong deformatin insertion actually) followed by a
  fibration

  \begin{center}
    \includegraphics[trim = {38 0 0 30}]{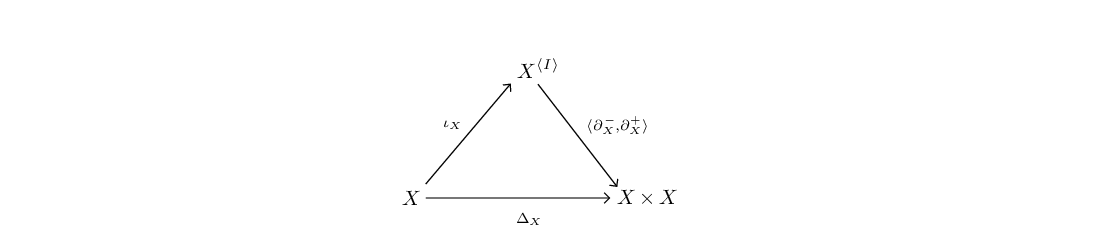}
  \end{center}

\end{remark}

\begin{definition}[Brown, 1973]
  A category $\mathbb{F}$ with finite limits equipped with a class of
  {\tmem{fibrations}} $\mathcal{F}$ and a class of {\tmem{weak equivalences}}
  $\mathcal{W}$ is a {\tmem{category of fibrant objects}} provided
  \begin{enumerateroman}
    \item $\tmop{Iso} (\mathbb{F}) \subset \mathcal{F} \cap \mathcal{W}$;

    \item $\mathcal{W}$ verifies 3-for-2;

    \item $\mathcal{F}$ and $\mathcal{F} \cap \mathcal{W}$ are closed under
    pullbacks;

    \item Any object is fibrant;

    \item for any $X \in \mathbb{F}$ there is a an object $\pao{X}$ such that
    there is a factorisation

    \begin{center}
      \includegraphics[trim = {38 0 0 30}]{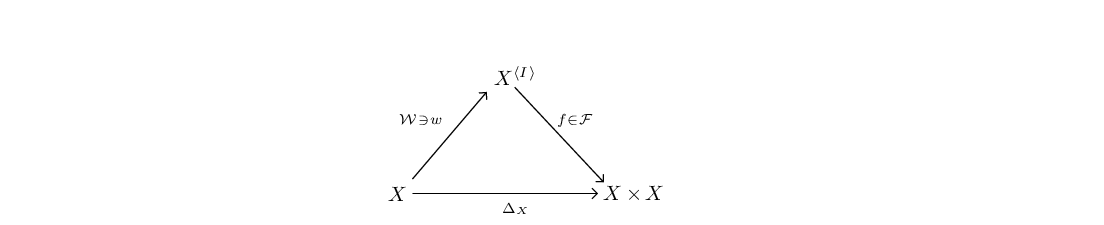}
    \end{center}
  \end{enumerateroman}
\end{definition}

\begin{theorem}
  Let $I$ be a Hurewicz interval. $(\mathbb{H}, I)$ with
  $\mathcal{F}=\mathcal{H}$ and $\mathcal{W}= \{ \tmop{homotopy}
  \tmop{equivalences} \}$ is a category of fibrant objects.
\end{theorem}

\begin{proof}

  \begin{enumerateroman}
    \item Obvious;

    \item Remark \ref{rem:three-for-two};

    \item Corollary \ref{cor:pullback} and Proposition
    \ref{prop:triv-pullback};

    \item Proposition \ref{prop:fibrant};

    \item Remark \ref{rem:diagonal-facto}.
  \end{enumerateroman}
\end{proof}

\section*{Appendix: the effective topos}

\begin{definition}
  An {\tmem{effective set}} $\eset{X}$ is a set equipped with an
  {\tmem{effective equality}}, that is a non-standard predicate
  $| - \approx -| : X \times X \rightarrow \mathcal{P}\mathbb{N}$ which is
  \begin{enumerateroman}
    \item {\tmem{symmetric:}} $\eeq{x}{x'}  \ent  \eeq{x'}{x}$

    \item {\tmem{transitive:}} $\eeq{x}{x'} \wedge \text{ }
    \eeq{x'}{x''}  \ent  \eeq{x}{x''}$
  \end{enumerateroman}
\end{definition}

\begin{remark}
  Notice that we do not alway have reflexivity, that is $\isval{}
  \eeq{x}{x}$. In fact, the latter assert's $x$'s {\tmem{existence}}.
  Accordingly, $E (x) \deq \eeq{x}{x} $ is called the {\tmem{existence
  predicate}} on $\eset{X}$. Call $x \in X$ {\tmem{ghost}} if its existence is
  empty. In particular, two equal inhabitants cannot be ghosts as
  \[ \eeq{x}{x'}  \ent E (x) \wedge E (x') \]
\end{remark}

\begin{definition}
  Let $\eset{X}$ and $\eset{Y}$ be effective sets. A \ functional relation
  \[ \Phi : \eset{X} \rightsquigarrow \eset{Y} \]
  is a predicate $\Phi : X \times Y \rightarrow \pn$ which is
  \begin{enumerateroman}
    \item {\tmem{extensional:}} $\Phi (x, y) \wedge \eeq{x}{x'} \wedge
    \eeq{y}{y'}  \ent \Phi (x', y')$

    \item {\tmem{strict:}} $\Phi (x, y)  \ent E (x) \wedge E (y)$

    \item {\tmem{single-valued:}} $\Phi (x, y) \wedge \Phi (x, y')  \ent
    \eeq{y}{y'}$

    \item {\tmem{total:}} $E (x)  \ent  \bigcup_{y \in Y} E (y) \wedge \Phi
    (x, y)$
  \end{enumerateroman}
  Two functional relations $\Phi, \Psi : \eset{X} \rightsquigarrow \eset{Y}$
  are equivalent if $\Phi \dashv \ent \Psi$.
\end{definition}

\begin{definition}
  A {\tmem{morphism}} of effective sets is an equivalence class of functional
  relations.
\end{definition}

\begin{notation}{We shall write $\rep{f}$ for an arbitrary but fixed representant
of the morphism $f$. \ }
\end{notation}

\begin{theorem}[Hyland]
  Effective sets and their morphisms aggregate to the category $\eff$ where
  \begin{enumeratenumeric}
    \item the composition of $f : (X, \approx) \rightarrow (Y, \approx)$ and
    $g : (Y, \approx) \rightarrow (Z, \approx)$ is represented by
    \begin{eqnarray*}
      \left( \rep{g} \circ \rep{f} \right) (x, z) & \deq & \bigcup_{y \in Y} E
      (y) \wedge \rep{f} (x, y) \wedge \rep{g} (y, z)
    \end{eqnarray*}
    \item the identity $\tmop{id}_{(X, \approx)}$ is represented by
    \[ \rep{\tmop{id}_X} (x, x') \assign \eeq{x}{x'} \]
  \end{enumeratenumeric}
  This category is a topos.
\end{theorem}

\begin{remark}
  Composition is $\exists y. \rep{f} (x, y) \wedge \rep{g} (y, z)$ in $\eff$'s
  internal logic.
\end{remark}

\begin{remark}
  \label{rem:Delta}Assume $(X, \approx) \in \eff$. The relation $\approx$ is
  an equivalence relation on
  \[ \underline{X}  \deq  \{ x \in X|E (x) \neq \varnothing \} \]
  Call {\tmem{equality class}} an element of the quotient $\underline{X} /
  \approx$. The assignment
  \begin{eqnarray*}
    \Gamma : (X, \approx) & \mapsto & \underline{X} / \approx
  \end{eqnarray*}
  extends to a functor $\Gamma : \eff \rightarrow \text{{\tmstrong{Set}}}$.
  This functor is isomorphic to the standard {\tmem{global sections functor}}
  and has a right adjoint $\Delta$, given on objects by
  \begin{eqnarray*}
    \Delta : A & \mapsto & (A, \approx_{\Delta})
  \end{eqnarray*}
  where $\approx_{\Delta}$ is the {\tmem{non-standard equality}}
  \begin{eqnarray*}
    a \approx_{\Delta} a' & = & \left\{\begin{array}{lll}
      \mathbb{N} &  & a = a'\\
      \varnothing &  & \tmop{otherwise}
    \end{array}\right.
  \end{eqnarray*}
  The pair $(\Gamma, \Delta)$ is a geometric morphism.
\end{remark}

\begin{definition}
  An {\tmem{assembly}} or {\tmem{$\omega$-set}} is an object $\eset{X} \in
  \eff$ such that $\eeq{x}{x'} = \varnothing$ if $x \neq x'$.
\end{definition}

\begin{remark}
  An assembly $(X, \approx)$ is thus given by the datum $(X, E)$ where $E : X
  \rightarrow \mathcal{P}\mathbb{N}$ is the existence predicate.
\end{remark}

\begin{example}
  {\tmdummy}

  \begin{enumeratenumeric}
    \item The assembly $\left( \left\{ \astt \right\}, \mathbb{N} \right)$ is
    terminal in $\eff$.

    \item The assembly $(\mathbb{N}, E(n) = \{ n \})$ is an NNO in {\eff}.

    \item Non-example: $\Omega$.
  \end{enumeratenumeric}
\end{example}

\begin{definition}
  Let $\eset{X}, \eset{Y} \in \eff$. A function $f : X \rightarrow Y$ is
  {\tmem{effective}} is there is a {\tmem{tracker}} $t \in \mathbb{N}$ such
  that , for all $x, x' \in X$ and $n \in \eeq{x}{x'}$ we have $\app{t}{n \in
  \eeq{f (x)}{f (x')}}$.
\end{definition}

\begin{remark}
  A effective function induces a morphism $\eset{X} \rightarrow \eset{Y}$
  represented by
  \[ \mathfrak{R}_f (x, y) = \bigcup_{x' \in X} \left\{ \langle m, n \rangle |
     m \in \eeq{x}{x'} \nocomma, n \in \eeq{f (x')}{y} \right\} \]
\end{remark}

\begin{proposition}
  Any morphism to an assembly is induced by a unique effective function.
\end{proposition}

\begin{corollary}
  A morphism $f : \eset{X} \rightarrow \eset{Y}$ among assemblies is induced
  by a {\tmem{supereffective}} function $f : X \rightarrow Y$ for which there
  exists a {\tmem{tracker}} $t \in \mathbb{N}$ such that we have $\app{t}{n
  \in E (f (x))}$ for all $x \in X$ and $n \in E (x)$.
\end{corollary}

\end{document}